\documentclass{amsart}
\copyrightinfo{2006}{American Mathematical Society}
\usepackage{mathrsfs,amscd}
\usepackage{amssymb}
\usepackage{amssymb, amsmath, mathrsfs}
\usepackage{amscd}
\usepackage{latexsym}
\usepackage{amssymb}
\usepackage{amsmath}
\usepackage{amsthm}
\usepackage{indentfirst}

\newtheorem{theorem}{Theorem}[section]

\theoremstyle{definition}

\theoremstyle{remark}

\numberwithin{equation}{section}
\allowdisplaybreaks%%%%%%%%%?????????????
\begin{document}

\title{Explicit singular minimal surface solutions for gravitational instantons}

\author[W. Yan]{Weiping Yan}
\address{School of Mathematics, Xiamen University, Xiamen 361000, People's Republic of China}
\address{Laboratoire Jacques-Louis Lions, Université Pierre et Marie Curie (Sorbonne Université), 4, Place Jussieu, 75252 Paris, France.}\email{yanwp@xmu.edu.cn}
\thanks{The author is supported by NSFC No 11771359.}

%\subjclass[2010]{74H35,35L05,35C06}
\date{Dec 2017}
%\keywords{Born-Infeld equation, Membrane equation, Self-s}

\begin{abstract}
We construct a family of instanton metric obtained from new exact singular solutions for minimal surfaces by noticing the correspondence between minimal surfaces in the three dimesional Euclidean space and gravitational instantons possessing two killing vectors. 
By Calabi's correspondence, we derive a family of explicit maximal surface solution for spacelike surface with zero mean curvature equation.
\end{abstract}
\maketitle

%\tableofcontents

\section{Introduction}
Gravitational instantons can be seen as localized ``excitations'' in imaginary time and mediate the transition between two distinct gravitational vacua \cite{Nu2,Ati,Hit,Eg,Chen}. It plays an important role in the Euclidean path integral approach to quantum gravity \cite{Haw,Gi1,Gi2,Gi3}. On one hand, an instanton solutions of the Einstein field equations with Euclidean signature and anti-self-dual curvature can obtain from the minimal surface in Euclidean $E^3(x,y,z)$ \cite{Nu1}.  
On the other hand, from J$\ddot{o}$rgens' correspondence \cite{J},
any two dimensional minimal surfaces gives a solution to the real elliptic Monge-Amp$\grave{e}$re equation on a real manifold of two dimensions. These solutions then lead to some gravitational instantons, which are described by hyper- K$\ddot{a}$hler metrics were studied extensively in the framework of supergravity and M-theory as well as Seiberg-Witten theory \cite{Gi4,SW}.

The minimal surface equation is the equation of minimal graphs over a domain of the $xy$-plane in $E^3$, which
takes the form
\begin{equation}\label{ENNN1-1}
\partial_x\left(\frac{\partial_x u}{\sqrt{1+|\partial_x u|^2+|\partial_yu|^2}}\right)+\partial_y\left(\frac{\partial_y u}{\sqrt{1+|\partial_x u|^2+|\partial_yu|^2}}\right)=0,
\end{equation}
it is equivalent to the classcial form 
\begin{equation*}
u_{xx}(1+u_y^2)-2u_xu_yu_{xy}+u_{yy}(1+u_x^2)=0,
\end{equation*}
which exhibits the following scaling invariance for any $\lambda>0$,
\begin{equation*}
u(x,y)\mapsto u_{\lambda}(x,y)=\lambda^{-1} u(\lambda x,\lambda y).
\end{equation*}
A class of gravitational instantons may be represented by the metric
\begin{equation}\label{E1-3}
ds^2=\frac{1+u_t^2}{\sqrt{1+u_t^2+u_x^2}}(dt^2+dy^2)+\frac{1+u_x^2}{\sqrt{1+u_t^2+u_x^2}}(dx^2+dz^2)+\frac{2u_tu_x}{\sqrt{1+u_t^2+u_x^2}}(dtdx+dydz),
\end{equation}
for which the Einstein field equations reduce to
\begin{equation}\label{E1-4}
u_{tt}(1+u_x^2)-2u_tu_xu_{tx}+u_{xx}(1+u_t^2)=0.
\end{equation}
\begin{theorem}
The equation (\ref{E1-4}) has a family of exact singular solutions
\begin{equation*}
u(t,x)=k\ln|\frac{x}{t}+(1+\frac{x^2}{t^2})^{\frac{1}{2}}|,~~\forall k\in \mathbb{R}\setminus \{0\}.
\end{equation*}
Moreover,  there is a class of gravitational instantons can be represented by the metric
\begin{equation*}
ds^2=a(t,x)(dt^2+dy^2)+b(t,x)(dx^2+dz^2)+c(t,x)(dtdx+dydz),
\end{equation*}
where
\begin{eqnarray*}
a(t,x)&=&\frac{1+u_t^2}{\sqrt{1+u_t^2+u_x^2}}\\
&=&[1+\frac{k^2x^2((t^2+x^2)^{\frac{1}{2}}+x)^2}{t^2(t^2+x^2)((t^2+x^2)^{\frac{1}{2}}+x)^2}][1+\frac{k^2x^2((t^2+x^2)^{\frac{1}{2}}+x)^2}{t^2(t^2+x^2)((t^2+x^2)^{\frac{1}{2}}+x)^2}+\frac{k^2((t^2+x^2)^{\frac{1}{2}}+x)^2}{(x(t^2+x^2)^{\frac{1}{2}}+t^2)^2}]^{-\frac{1}{2}},\\
b(t,x)&=&\frac{1+u_x^2}{\sqrt{1+u_t^2+u_x^2}}\\
&=&[1+\frac{k^2((t^2+x^2)^{\frac{1}{2}}+x)^2}{(x(t^2+x^2)^{\frac{1}{2}}+t^2)^2}][1+\frac{k^2x^2((t^2+x^2)^{\frac{1}{2}}+x)^2}{t^2(t^2+x^2)((t^2+x^2)^{\frac{1}{2}}+x)^2}+\frac{k^2((t^2+x^2)^{\frac{1}{2}}+x)^2}{(x(t^2+x^2)^{\frac{1}{2}}+t^2)^2}]^{-\frac{1}{2}},\\
c(t,x)&=&\frac{2u_tu_x}{\sqrt{1+u_t^2+u_x^2}}\\
&=&-2k^2(x(t^2+x^2)^{\frac{1}{2}}+x^2)((t^2+x^2)^{\frac{1}{2}}+x)(tx(t^2+x^2)^{\frac{1}{2}}+t^3+tx^2)^{-1}(x(t^2+x^2)^{\frac{1}{2}}+t^2)^{-1}\\
&&\times[1+\frac{k^2x^2((t^2+x^2)^{\frac{1}{2}}+x)^2}{t^2(t^2+x^2)((t^2+x^2)^{\frac{1}{2}}+x)^2}+\frac{k^2((t^2+x^2)^{\frac{1}{2}}+x)^2}{(x(t^2+x^2)^{\frac{1}{2}}+t^2)^2}]^{-\frac{1}{2}}.
\end{eqnarray*}
\end{theorem}
Obviously, equation (\ref{E1-4}) is the minimal surface equation (\ref{ENNN1-1}). It is well-known that Nutku\cite{Nu1} gave a famility of singularity solution of minimal surface equation (\ref{E1-4}) as
\begin{equation*}
u(t,x)=k\arctan\frac{x}{t},~~k\in\mathbb{R}\setminus \{0\}.
\end{equation*}

In this note, we show a new family of singular solutions of minimal surface equation (\ref{E1-4}), which gives a new of gravitational instantons.

\section{Exact gravitational instantons}
Introduce the similarity coordinates
\begin{equation*}
\tau=-\log t,~~~~\rho=\frac{x}{t},
\end{equation*}
then we denote by
\begin{equation*}
u(t,x)=v(-\log t,\frac{x}{t}),
\end{equation*}
and noticing 
\begin{eqnarray*}
&&u_t(t,x)=e^{\tau}(v_{\tau}+\rho v_{\rho}),\\
&&u_{tt}(t,x)=e^{2\tau}(v_{\tau\tau}+v_{\tau}+2\rho v_{\rho}+2\rho v_{\tau\rho}+\rho^2v_{\rho\rho}),\\
&&u_x(t,x)=e^{\tau}v_{\rho},\\
&&u_{xx}(t,x)=e^{2\tau}v_{\rho\rho},\\
&&u_{tx}(t,x)=e^{2\tau}(v_{\tau\rho}+v_{\rho}+\rho v_{\rho\rho})
\end{eqnarray*}
equation (\ref{E1-4}) is transformed into an one dimensional quasilinear elliptic equation
\begin{eqnarray}\label{YNN1-1}
v_{\tau\tau}+(1+\rho^2)v_{\rho\rho}&+&v_{\tau}+2\rho v_{\rho}+2\rho v_{\tau\rho}+e^{2\tau}v_{\rho}^2(v_{\tau\tau}+v_{\tau}+2\rho v_{\rho}+2\rho v_{\tau\rho}+\rho^2v_{\rho\rho})\nonumber\\
&&+e^{2\tau}(v_{\tau}+\rho v_{\rho})^2v_{\rho\rho}-2e^{2\tau}v_{\rho}(v_{\tau}+\rho v_{\rho})(v_{\rho}+\rho v_{\rho\rho}+v_{\tau\rho})=0.
\end{eqnarray}
The equation only on $\rho$ of elliptic equation (\ref{YNN1-1}) is 
\begin{equation*}
(\rho^2+1)v_{\rho\rho}+2\rho v_{\rho}=0,
\end{equation*}
which is an ODE. Direct computation shows that it has a family of solutions
\begin{equation*}
v(\rho)=k\ln|\rho+(1+\rho^2)^{\frac{1}{2}}|,
\end{equation*}
where $k$ is an arbitrary constant in $\mathbb{R}\setminus \{0\}$.

Minimal surface equation (\ref{E1-4}) has a family of explicit self-similar solutions
\begin{equation*}
u(t,x)=k\ln|\frac{x}{t}+(1+\frac{x^2}{t^2})^{\frac{1}{2}}|,~~\forall k\in \mathbb{R}\setminus \{0\}.
\end{equation*}
It is easy to see that
\begin{equation*}
\partial_xu|_{x=0}=\frac{k}{t}\rightarrow+\infty,~~as~~t\rightarrow 0.
\end{equation*}
So by (\ref{E1-3}), a class of gravitational instantons can be represented by the metric
\begin{equation*}
ds^2=a(t,x)(dt^2+dy^2)+b(t,x)(dx^2+dz^2)+c(t,x)(dtdx+dydz),
\end{equation*}
where
\begin{eqnarray*}
a(t,x)&=&\frac{1+u_t^2}{\sqrt{1+u_t^2+u_x^2}}\\
&=&[1+\frac{k^2x^2((t^2+x^2)^{\frac{1}{2}}+x)^2}{t^2(t^2+x^2)((t^2+x^2)^{\frac{1}{2}}+x)^2}][1+\frac{k^2x^2((t^2+x^2)^{\frac{1}{2}}+x)^2}{t^2(t^2+x^2)((t^2+x^2)^{\frac{1}{2}}+x)^2}+\frac{k^2((t^2+x^2)^{\frac{1}{2}}+x)^2}{(x(t^2+x^2)^{\frac{1}{2}}+t^2)^2}]^{-\frac{1}{2}},\\
b(t,x)&=&\frac{1+u_x^2}{\sqrt{1+u_t^2+u_x^2}}\\
&=&[1+\frac{k^2((t^2+x^2)^{\frac{1}{2}}+x)^2}{(x(t^2+x^2)^{\frac{1}{2}}+t^2)^2}][1+\frac{k^2x^2((t^2+x^2)^{\frac{1}{2}}+x)^2}{t^2(t^2+x^2)((t^2+x^2)^{\frac{1}{2}}+x)^2}+\frac{k^2((t^2+x^2)^{\frac{1}{2}}+x)^2}{(x(t^2+x^2)^{\frac{1}{2}}+t^2)^2}]^{-\frac{1}{2}},\\
c(t,x)&=&\frac{2u_tu_x}{\sqrt{1+u_t^2+u_x^2}}\\
&=&-2k^2(x(t^2+x^2)^{\frac{1}{2}}+x^2)((t^2+x^2)^{\frac{1}{2}}+x)(tx(t^2+x^2)^{\frac{1}{2}}+t^3+tx^2)^{-1}(x(t^2+x^2)^{\frac{1}{2}}+t^2)^{-1}\\
&&\times[1+\frac{k^2x^2((t^2+x^2)^{\frac{1}{2}}+x)^2}{t^2(t^2+x^2)((t^2+x^2)^{\frac{1}{2}}+x)^2}+\frac{k^2((t^2+x^2)^{\frac{1}{2}}+x)^2}{(x(t^2+x^2)^{\frac{1}{2}}+t^2)^2}]^{-\frac{1}{2}}.
\end{eqnarray*}

\section{Discussion}
In section 2, we derive a singularity metric, the singular point is at $t=0$. In fact, let $T$ be a positive parameter. if we introduce the similarity coordinates
\begin{equation*}
\tau=-\log(T-t),~~~~\rho=\frac{x}{T-t},
\end{equation*}
then minimal surface equation (\ref{E1-4}) has a family of explicit self-similar singularity solutions
\begin{equation*}
u(t,x)=k\ln|\frac{x}{T-t}+(1+\frac{x^2}{(T-t)^2})^{\frac{1}{2}}|,~~\forall k\in \mathbb{R}\setminus \{0\}.
\end{equation*}
It is easy to see that
\begin{equation*}
\partial_xu|_{x=0}=\frac{k}{T-t}\rightarrow+\infty,~~as~~t\rightarrow T^{-1}.
\end{equation*}
Hence we have a class of gravitational instantons can be represented by the metric
\begin{equation*}
ds^2=a(T-t,x)(dt^2+dy^2)+b(T-t,x)(dx^2+dz^2)+c(T-t,x)(dtdx+dydz).
\end{equation*}

On the other hand, using Calabi's correspondence, we know that spacelike surface with zero mean curvature equation
\begin{equation*}
u_{tt}(1-u_x^2)+2u_tu_xu_{tx}+u_{xx}(1-u_t^2)=0
\end{equation*}
has the same self-similar solutions with minimal surface equation (\ref{E1-4}).

%At last, we propose two question which will be considered in the furture.

%ii) Can we give new exact singularity solutions to Einstein-Maxwell theory from our exact singularity solutions of minimial surface?

%\textbf{Acknowledgments.} The author express his thanks to Prof.C. Xia for his discussion on minimal surface.

\end{document}